# INDIVIDUAL VERSUS CLUSTER RECOVERIES WITHIN A SPATIALLY STRUCTURED POPULATION

By L. Belhadji and N. Lanchier

*CNRS—Université de Rouen*

Stochastic modeling of disease dynamics has had a long tradition. Among the first epidemic models including a spatial structure in the form of local interactions is the contact process. In this article we investigate two extensions of the contact process describing the course of a single disease within a spatially structured human population distributed in social clusters. That is, each site of the $d$-dimensional integer lattice is occupied by a cluster of individuals; each individual can be healthy or infected. The evolution of the disease depends on three parameters, namely the outside infection rate which models the interactions between the clusters, the within infection rate which takes into account the repeated contacts between individuals in the same cluster, and the size of each social cluster. For the first model, we assume cluster recoveries, while individual recoveries are assumed for the second one. The aim is to investigate the existence of nontrivial stationary distributions for both processes depending on the value of each of the three parameters. Our results show that the probability of an epidemic strongly depends on the recovery mechanism.

**1. Introduction.** To study the course of a disease within a spatially structured population, Harris [8] introduced a model known as the basic contact process. Each site of the $d$-dimensional integer lattice is occupied by an individual either healthy or infected; each individual gets infected at a rate that depends on the number of infected individuals in some interaction neighborhood. Including spatial structure in the form of local interactions has shown that, for a disease to spread, the infection rate needs to exceed a threshold that is greater than the one for a nonspatial population. The smaller the size of the interaction neighborhood, the greater the threshold for the disease to spread. The reason for this is the lack of sufficient numbers of susceptible individuals near the location of a disease outbreak once the disease starts spreading.

The study of the evolution of diseases within spatially structured populations, including extensions of the contact process, is widespread in the particle system literature. The first process we investigate in this article has been introduced by Schinazi [12], and will be referred to as the *cluster recovery process* (CRP). His model studies the spread of an infectious disease such as tuberculosis within a population grouped in social clusters, each of the clusters having the same size. The second model we investigate assumes another recovery mechanism and will be referred to as the *individual recovery process* (IRP).

For both CRP and IRP, the dynamics depends on three parameters, namely the outside infection rate $\lambda$ (the rate at which an individual infects healthy individuals of other clusters), the within infection rate $\phi$ (the rate at which an individual infects healthy individuals present in the same cluster) and the cluster size $\kappa$ (which is constant regardless of the location of the cluster on the lattice). The population is divided into social clusters, the individuals belonging to the same cluster having repeated contacts and the individuals belonging to neighboring clusters having casual contacts only, which suggests that the infection spreads out faster *vertically* within the clusters than *horizontally* between the clusters. In particular, it is assumed for both models that, once a cluster has at least one infected individual, infections within the cluster are much more likely than additional infections from the outside so we neglect the latter. The cluster size can be seen as the mean number of individuals having sustained contacts with a given individual. Even if, in a more realistic setting, this parameter should fluctuate significantly depending on the social customs of each of the individuals, it is assumed for technical reasons that all the clusters have the same size. In particular, whereas the infection rates $\lambda$ and $\phi$ are parameters linked to the nature of the disease, the cluster size $\kappa$ depends on the social customs of the individuals. The only difference between the CRP and the IRP is the recovery mechanism. For the CRP, all the infected individuals in a given cluster simultaneously recover at rate 1 thanks to an antidote. This applies to places where there is a good tracking system of infectious diseases so that, once an infected individual is discovered, its social cluster is rapidly tracked down [12]. For the IRP, we deal with the other extreme case when at most one infected individual recovers at once; that is, the tracking system is not effective enough and the infection can spread within a given cluster before it is detected. In particular, the CRP and the IRP can be considered as spatial stochastic models for the transmission of infectious diseases in developed and developing countries, respectively. As we will see further, this difference implies that, the cluster size being fixed, an epidemic may occur for the IRP provided the outside infection rate is strictly positive. In the CRP, the whole population always recovers if the outside infection rate is smaller than some critical value depending on the cluster size. Moreover, we obtain for the IRP that, when the within infection rate is greater than 1, even if the outside infection rate is low, an epidemic is possible provided the cluster size is large enough. The CRP exhibits the opposite behavior in the sense that, the within infection rate being fixed, if the outside infection rate is low, there is no epidemic whatever the cluster size.

To figure out the differences between individual and cluster recoveries, we now introduce the explicit dynamics of both processes. The IRP is a continuous-time Markov process in which the state at time $t$ is a function $\xi_t : \mathbb{Z}^d \longrightarrow \{0, 1, \ldots, \kappa\}$, with $\kappa$ denoting the common size of the clusters. The cluster at site $x \in \mathbb{Z}^d$ is said to be healthy at time $t \geq 0$ if $\xi_t(x) = 0$, and infected otherwise. More precisely, $\xi_t(x)$ indicates the number of infected individuals present in the cluster at time $t \geq 0$. To take into account the outside infections, we also introduce an interaction neighborhood. For any $x, z \in \mathbb{Z}^d$, $x \sim z$ indicates that site $z$ is one of the $2d$ nearest neighbors of site $x$. Then, the state of site $x$ flips according to the following transition rates:

$$
\begin{aligned}
0 \to 1 \quad &\text{at rate} \quad \lambda \sum_{x \sim z} \xi(z), \\
i \to i+1 \quad &\text{at rate} \quad i\phi, \quad i = 1, 2, \ldots, \kappa - 1, \\
i \to i-1 \quad &\text{at rate} \quad i, \quad i = 1, 2, \ldots, \kappa.
\end{aligned}
$$

In other words, a healthy cluster at site $x$ gets infected, that is, the state of $x$ flips from 0 to 1, at rate $\lambda$ times the number of infected individuals present in the neighboring clusters. In other respects, if there are $i$ infected individuals in the cluster, $i = 1, 2, \ldots, \kappa - 1$, then each of these individuals infects healthy individuals in the cluster at rate $\phi$. Finally, each infected individual recovers at rate 1 regardless of the number of infected individuals in its cluster.

The CRP is a Markov process $\eta_t : \mathbb{Z}^d \longrightarrow \{0, 1, \ldots, \kappa\}$, with $\eta_t(x)$ denoting the number of infected individuals at site $x$ at time $t \geq 0$, and whose dynamics is obtained by replacing the transitions $i \to i-1$, $i = 1, 2, \ldots, \kappa$, above by the transitions

$$i \to 0 \quad \text{at rate} \quad 1, \quad i = 1, 2, \ldots, \kappa.$$

That is, all the infected individuals in a given cluster are now simultaneously replaced by healthy ones at rate 1, the infection mechanism being unchanged. For more details, see [12].

*The graphical representation.* An argument of Harris [7] assures us of the existence and uniqueness of our spatially explicit, stochastic models. The idea is to construct the processes from collections of independent Poisson processes, which is referred to as Harris' graphical representation. For each $x$, $z \in \mathbb{Z}^d$ with $x \sim z$ and $i = 1, 2, \ldots, \kappa$, we let $\{T_n^{x,z,i} : n \geq 1\}$ denote the arrival times of independent Poisson processes with rate $\lambda$, and draw an arrow labeled with an $i$ from site $x$ to site $z$ at time $T_n^{x,z,i}$ to indicate that an outside infection may occur. To take into account the within infections, we introduce, for $x \in \mathbb{Z}^d$ and $i = 1, 2, \ldots, \kappa - 1$, a further collection of independent Poisson processes, denoted by $\{U_n^{x,i} : n \geq 1\}$; each of them has rate $\phi$. We put the symbol $\bullet_i$ at $(x, U_n^{x,i})$ to indicate that an infection from the inside may occur. Finally, for each $x \in \mathbb{Z}^d$ and $i = 1, 2, \ldots, \kappa$, we let $\{V_n^{x,i} : n \geq 1\}$ be the arrival times of independent rate-1 Poisson processes, and put a $\times_i$ at site $x$ to indicate that a recovery may occur.

Given initial configurations $\xi_0$ and $\eta_0$, and the graphical representation introduced above, the processes can be constructed as follows. If there are at least $i$ infected individuals at site $x$ at time $T_n^{x,z,i}$ and the cluster at $z$ is healthy, then the state of $z$ flips from 0 to 1 for both processes. In other respects, if there are $j$ infected individuals, $i \leq j \leq \kappa - 1$, at site $x$ at time $U_n^{x,i}$, then one more individual gets infected in the cluster, that is, the state of $x$ flips from $j$ to $j + 1$, for both processes. Finally, if there are $j$ infected individuals, $1 \leq j \leq \kappa$, at site $x$ at time $V_n^{x,i}$, then the state of $x$ flips from $j$ to $j - 1$ if and only if $i \leq j$ for the IRP, while it flips from $j$ to 0 if and only if $i = 1$ for the CRP. In particular, $\times_i$'s, $i = 2, 3, \ldots, \kappa$, have no effect on the CRP.

A nice feature of the graphical representation is that it allows us to couple several processes starting from different initial configurations, which can be done by using the same collections of Poisson processes. See [5], page 119 and [10], page 32.

*The mean-field model of the IRP.* To figure out the properties of our spatial models, the first step is to investigate their deterministic nonspatial versions called the mean-field model; that is, we assume that all sites are independent and the system is spatially homogeneous. This then results in a system of ordinary differential equations for the densities of healthy and infected clusters. The reason for introducing this model is that the existence of locally stable fixed points for the mean-field model may be symptomatic of the existence of stationary measures for the original spatial model. See [6] for different possible studies of a model. Let $u_i$ denote the density of clusters with $i$ infected individuals, $i = 0, 1, \ldots, \kappa$. The mean-field model of the IRP is then described by the following coupled system of ordinary differential equations:

$$\frac{du_0}{dt} = u_1 - \lambda u_0 \sum_{i=1}^{\kappa} i u_i,$$

$$\frac{du_1}{dt} = \lambda u_0 \sum_{i=1}^{\kappa} i u_i - (1+\phi)u_1 + 2u_2,$$

$$\frac{du_i}{dt} = (i-1)\phi u_{i-1} - i(1+\phi)u_i + (i+1)u_{i+1}, \qquad i=2,3,\ldots,\kappa-1,$$

$$\frac{du_\kappa}{dt} = (\kappa-1)\phi u_{\kappa-1} - \kappa u_\kappa.$$

To find the condition for the existence of a nontrivial equilibrium, we start by setting the right-hand side of each of the equations equal to 0. From the last equation, we obtain

$$u_\kappa = \frac{\kappa-1}{\kappa}\phi u_{\kappa-1}.$$

By induction, we get

$$u_i = \frac{i-1}{i}\phi u_{i-1}, \qquad i=2,3,\ldots,\kappa,$$

from which it follows that

$$u_i = \frac{\phi^{i-1}}{i} u_1, \qquad i=2,3,\ldots,\kappa.$$

Reporting in the first equation then leads to

$$\frac{du_0}{dt} = u_1 - \lambda \sum_{i=0}^{\kappa-1} \phi^i u_0 u_1 = 0.$$

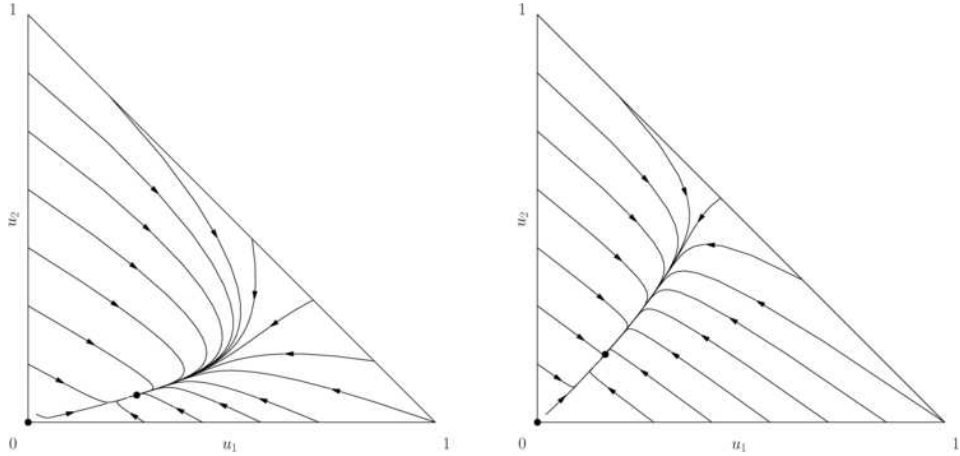

FIG. 1. *Solution curves of the mean-field model of the IRP in the case $\kappa=2$. Left: $\lambda=1$ and $\phi=0.5$. Right: $\lambda=0.5$ and $\phi=2$.*

The previous equation tells us that the condition for an epidemic to occur is given by

$$\lambda(1 + \phi + \cdots + \phi^{\kappa-1}) > 1. \tag{1}$$

When (1) does not hold, the boundary equilibrium $u_0 = 1$ and $u_i = 0$ for $i = 1, 2, \ldots, \kappa$ is the unique fixed point of the system of ordinary differential equations. When (1) is satisfied, there is a nontrivial equilibrium characterized by

$$u_0 = \frac{1}{\lambda} \times \left\{\sum_{i=0}^{\kappa-1} \phi^i\right\}^{-1} \quad \text{and} \quad u_i = \frac{\phi^{i-1}}{i} u_1, \qquad i = 2, 3, \ldots, \kappa.$$

First of all, note that, when $\lambda > 1$, condition (1) holds for all the values of $\phi \geq 0$ and $\kappa \geq 1$. In other respects, when $\lambda \leq 1$ and $\kappa \geq 2$, there is a nontrivial equilibrium provided $\phi$ is sufficiently large. In the same way, when $\lambda \leq 1$ and $\phi \geq 1$, there is a nontrivial equilibrium if the cluster size $\kappa$ is large enough. Moreover, numerical simulations indicate that, when condition (1) holds, the nontrivial equilibrium is locally stable, which suggests the existence of a nontrivial stationary distribution for the corresponding spatial model. See Figure 1 for a picture. For a description of the mean-field model of the CRP, see [12], Section 3.

*The stochastic models.* We now discuss the effects of each of the three parameters, namely the outside infection rate $\lambda$, the cluster size $\kappa$ and the within infection rate $\phi$, on the probability of an epidemic for both models. We will prove that the behavior predicted by the mean-field model above holds as well for the IRP, and provide comparisons between both processes. Except for Theorem 6 whose proof is carried out in this paper, all the results regarding the CRP can be found in [12], so we now focus especially on the IRP, but the following definitions and basic properties hold for the CRP as well.

We say that an epidemic may occur for the IRP with parameters $(\lambda, \kappa, \phi)$ if

$$P_\xi^0(|\xi_t| \geq 1 \text{ for all } t \geq 0) > 0,$$

where $P_\xi^0$ denotes the law of the process starting from one infected individual at the origin, and where $|\xi_t|$ is the total number of infected individuals at time $t \geq 0$. The previous property is equivalent to the existence of a nontrivial stationary distribution, where nontrivial means $\neq \delta_0$, the "all 0" configuration. We say, on the contrary, that there is no epidemic if $\delta_0$ is the unique stationary distribution.

As in [12], one of the keys of our results is monotonicity. A crucial feature of the CRP that ensures monotonicity is the fact that the transition $i \to 0$ occurs at a constant rate (actually, any transition rate which is decreasing in $i$ would work). For the IRP, this follows from the fact that we have at most one recovery at any fixed time. The following monotonicity result can be checked by using a coupling argument if we think of the process as being generated by Harris' graphical representation introduced above.

LEMMA 1.1. *The IRP is attractive and monotone with respect to $\lambda$, $\kappa$ and $\phi$.*

Here attractivity means that if $\xi_0^1(x) \leq \xi_0^2(x)$ for any $x \in \mathbb{Z}^d$ at time 0, then $\xi_t^1$ and $\xi_t^2$ can be constructed in the same probability space in such a way that

$$P^{1,2}(\forall\, x \in \mathbb{Z}^d, \xi_t^1(x) \leq \xi_t^2(x)) = 1 \qquad \text{for any } t \geq 0,$$

with $P^{1,2}$ denoting the law of the coupled process starting from $(\xi_0^1, \xi_0^2)$.

In order to introduce our results, we first observe that, when $\kappa = 1$, each of the clusters is only in one of the two states $0 =$ healthy or $1 =$ infected, and the process $\xi_t$ reduces to a basic contact process with parameter $\lambda$. In this case, there exists a critical value $\lambda_c \in (0, \infty)$ such that if $\lambda \leq \lambda_c$, then the process converges in distribution to the "all 0" configuration; otherwise, an epidemic may occur. See [2] or [10], Theorem 2.25. This, together with Lemma 1.1, implies that, when $\lambda > \lambda_c$, an epidemic may occur for any $\kappa \geq 1$ and $\phi \geq 0$. This case corresponds to the case $\lambda > 1$ of the mean-field model.

To figure out the intuition behind our first result, we start by removing the interactions between clusters by setting $\lambda = 0$. This makes the IRP a system of independent random walks with absorbing state 0; each of them represents the number of infected individuals in the associated cluster. Since each of these random walks returns to 0 with probability 1, the process converges in distribution to the "all 0" configuration. By relying on a perturbation argument, the result can be extended to the region $\lambda > 0$ small. More precisely, we have the following

THEOREM 1. *For all $\kappa \geq 1$ and $\phi \geq 0$, there is $\lambda_c(\kappa, \phi) \in (0, \infty)$ such that if $\lambda < \lambda_c(\kappa, \phi)$ there is no epidemic for the IRP, while if $\lambda > \lambda_c(\kappa, \phi)$ an epidemic may occur.*

Note that, due to the monotonicity with respect to $\kappa$, we get $\lambda_c(\kappa, \phi) \leq \lambda_c(1, \phi) = \lambda_c$. The analogue of Theorem 1 for the CRP is given by

THEOREM 2 ([12]). *For all $\kappa \geq 1$ and $\phi \geq 0$, if $\lambda \leq \lambda_c/\kappa$ there is no epidemic for the CRP, while if $\lambda > \lambda_c$ an epidemic may occur.*

The next step is to set $\lambda < \lambda_c$ and $\kappa \geq 2$, and to discuss the probability of an epidemic depending on the value of the within infection rate $\phi$. To begin with, when $\phi = 0$ the IRP reduces to a basic contact process with parameter $\lambda$. In particular, since $\lambda < \lambda_c$, there is no epidemic. The other extreme formal case $\phi = \infty$ corresponds to a Richardson model with parameter $\lambda\kappa$ so that an epidemic may occur provided $\lambda > 0$. To see that the previous two conclusions still hold when $\phi \in (0, \infty)$, we will first rely on a rescaling argument to estimate the rate of convergence of $P_\xi(\xi_t(x) = 0)$ in the two limiting cases $\phi = 0$ and $\phi = \infty$, respectively, where $P_\xi$ denotes the law of the IRP. These estimates will have to be good enough so that a perturbation argument can be applied. This, together with Lemma 1.1, implies that

THEOREM 3.   *For all $\kappa \geq 2$ and $\lambda \in (0, \lambda_c)$, there is $\phi_c(\lambda, \kappa) \in (0, \infty)$ such that if $\phi < \phi_c(\lambda, \kappa)$ there is no epidemic for the IRP, while if $\phi > \phi_c(\lambda, \kappa)$ an epidemic may occur.*

The CRP exhibits a quite different behavior since such a critical value exists if and only if the outside infection rate belongs to $(\lambda_c/\kappa, \lambda_c)$, that is:

THEOREM 4 ([12]).   *For all $\kappa \geq 2$ and $\lambda \in (\lambda_c/\kappa, \lambda_c)$, there is $\phi_c(\lambda, \kappa) \in (0, \infty)$ such that if $\phi < \phi_c(\lambda, \kappa)$ there is no epidemic for the CRP, while if $\phi > \phi_c(\lambda, \kappa)$ an epidemic may occur.*

The last step is to investigate the effects of the cluster size $\kappa$ on the probability of an epidemic for both models. For the mean-field model of the IRP, we have seen that, in the case $\lambda > 0$ and $\phi \geq 1$, condition (1) holds provided $\kappa$ is sufficiently large. The assumption $\phi \geq 1$ is to make sure that the sequence $1 + \phi + \cdots + \phi^{\kappa-1}$ diverges as $\kappa \to \infty$. The analogous result for the spatial version is given by the following

THEOREM 5.   *For all $\phi > 1$ and $\lambda \in (0, \lambda_c)$, there is $\kappa_c(\lambda, \phi) \geq 2$ such that if $\kappa < \kappa_c(\lambda, \phi)$ there is no epidemic for the IRP, while if $\kappa > \kappa_c(\lambda, \phi)$ an epidemic may occur.*

To understand the assumption $\phi > 1$, we observe that when the interactions between the clusters are removed, that is, when $\lambda = 0$, the process becomes a system of independent random walks with absorbing state 0. When $\phi > 1$ and only in this case, these random walks have a drift to the right so that the infection in a given cluster can persist a very long time. As we will see further, the previous observation is the key for proving Theorem 5. Our last result tells us that the CRP exhibits the opposite behavior in the sense that, $\phi$ being fixed, when $\lambda$ is too small, there is no epidemic whatever the cluster size. More precisely, we have the following

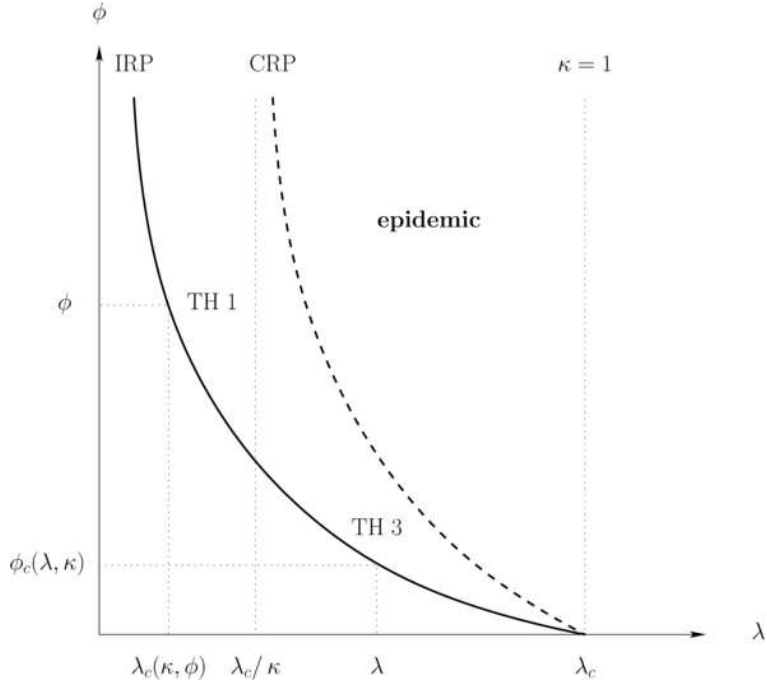

FIG. 2. *Phase diagram of the stochastic models. The continuous curve refers to the IRP and the dashed one to the CRP. For the IRP, Theorems* 1 *and* 3 *imply that the straight line* $\lambda = 0$ *is an asymptote of the phase transition curve, and Theorem* 5 *that, when* $\phi > 1$, $\lambda_c(\kappa, \phi) \to 0$ *as* $\kappa \to \infty$. *For the CRP, Theorem* 6 *implies that* $\lambda_c(\kappa, \phi) \to \lambda_c(\phi) > 0$ *as* $\kappa \to \infty$. *The monotonicity follows from Lemma* 1.1 *and Proposition* 1 *of* [12] *while we conjecture the continuity of the two curves.*

THEOREM 6. *For all $\phi \geq 0$ there exists $\lambda_c(\phi) \in (0, \infty)$ such that for any $\kappa \geq 1$ there is no epidemic for the CRP provided $\lambda \leq \lambda_c(\phi)$.*

The key idea of the proof is that, contrary to the IRP in which the number of infected individuals in a given cluster performs a random walk with a drift to the right when $\phi > 1$, for the CRP, the clusters globally recover sufficiently often so that the number of infected individuals between two recoveries cannot exceed some threshold $\kappa_0(\phi)$ with probability close to 1. In particular, by taking $\lambda > 0$ smaller than some critical value depending on $\kappa_0(\phi)$, but independent of $\kappa$, the disease dies out. All our results are summarized in Figure 2 above.

The rest of the article is devoted to proofs. In Section 2 we rely on a rescaling argument to prove Theorem 1. In Section 3 we use the fact that, when $\lambda \in (0, \lambda_c)$, the contact process with parameter $\lambda$ converges to the "all 0" configuration while the Richardson model with growth rate $\lambda \kappa$ has a nontrivial stationary distribution, to deduce Theorem 3. The proof of Theorem 5, partially based on random walks estimates, is carried out in Section 4. Finally, Section 5 is devoted to the proof of Theorem 6.

**2. Proof of Theorem 1.** In Section 1 we have seen that when $\lambda > \lambda_c$ an epidemic may occur whatever the size of the clusters and the value of the within infection rate. In view of the monotonicity of the IRP with respect to the outside infection rate $\lambda$, if in addition we prove that, $\kappa \geq 1$ and $\phi \geq 0$ being fixed, there exists $\lambda_0 > 0$ such that, when $\lambda \leq \lambda_0$, the "all 0" configuration is the unique stationary distribution, Theorem 1 will follow.

In order to prove the existence of such a $\lambda_0 > 0$, we use a rescaling argument to compare our stochastic process with an oriented percolation process on $\mathcal{G} = \mathbb{Z}^d \times \mathbb{Z}_+$. The rescaling argument has been invented by Bramson and Durrett [3] and is reviewed in [5]. To make the connection between the particle system and an oriented percolation process, the basic idea is to turn the space–time of the process into a brick wall, and to associate each brick with a certain good event $E_{z,n}$, $(z,n) \in \mathcal{G}$. For any $(z,n) \in \mathcal{G}$, we consider the vertical segment

$$S_{z,n} = \{(z,t) : nT \leq t \leq (n+1)T\} \tag{2}$$

as well as the space–time cylinder

$$C_{z,n} = \{(x,t) \in \mathbb{Z}^d \times \mathbb{R}_+ : x \sim z \text{ and } nT \leq t \leq (n+1)T\} \tag{3}$$

where $T$ is a large integer to be fixed later. We say that a site $(z,n) \in \mathcal{G}$ is good if the cluster at site $z$ is healthy at time $(n+1)T$. We want the good event $E_{z,n}$ to be measurable with respect to the graphical representation in $S_{z,n} \cup C_{z,n}$ and such that on $E_{z,n}$ site $(z,n)$ is good regardless of the value of $\xi_{nT}(z)$ or the configuration of the system in the cylinder $C_{z,n}$. We want this event to have a probability close to 1.

To define $E_{z,n}$ we start by setting $\lambda = 0$ to make our process a collection of independent random walks. Moreover, to be in the worst case, we assume that $\xi_{nT}(z) = \kappa$. Let $T_0^z$ be the first time $\xi_t(z)$ hits 0 after time $nT$. Since $\kappa < \infty$, the stopping time $T_0^z$ is finite with probability 1. This, together with the Beppo–Lévi theorem, implies that, for any $\varepsilon > 0$, there is a large enough $T$, fixed from now on, such that

$$P_\xi(T_0^z - nT \leq T) \geq 1 - \varepsilon/2, \tag{4}$$

where $P_\xi$ denotes the law of the particle system. To define the event $E_{z,n}$, we first require that $T_0^z \leq (n+1)T$ (which can be defined, in terms of Poisson processes, as a certain succession of ×'s and •'s at site $z$) so that on $E_{z,n}$ site $(z,n)$ is good provided $\lambda = 0$. To take into account the local interactions between the clusters when $\lambda > 0$, we complete the definition of $E_{z,n}$ by requiring no arrow to come from the cylinder $C_{z,n}$ to the segment $S_{z,n}$.

Finally, now that $T$ is fixed so that (4) holds, we can find $\lambda_0 > 0$ so small that

$$P_\xi((z,n) \text{ is good}) \geq P_\xi(E_{z,n}) \geq 1 - \varepsilon \qquad \text{if } \lambda \leq \lambda_0.$$

To complete the comparison, we now position oriented edges between sites in $\mathcal{G}$ in order to obtain an oriented percolation model. For $(x,n), (y,m) \in \mathcal{G}$, we draw an edge from $(x,n)$ to $(y,m)$ if and only if $x \sim y$ and $n = m$, or $x = y$ and $m = n + 1$. In that way, we may define a 1-dependent oriented percolation process on $\mathcal{G}$ for which sites are open with probability $1 - \varepsilon$. By choosing $\varepsilon > 0$ sufficiently small, the probability of a path of length $n$ of closed sites within this percolation process decreases exponentially fast with $n$ (see (8.2) in [1]). Since the existence of an infected site at time $(n+1)T$ implies the existence of a path of length $n$ of closed sites, Theorem 1 follows.

**3. Proof of Theorem 3.** In this section we first prove that, given $\lambda \in (0, \lambda_c)$ and $\kappa \geq 2$, for the spatial model with parameters $(\lambda, \kappa, \phi)$, there is no epidemic when $\phi = 0$ while an epidemic may occur when $\phi = \infty$. Then, by using a perturbation argument, we will deduce from the behavior of the IRP in the two extreme cases that there exist $\phi_1 > 0$ and $\phi_2 < \infty$ such that if $\phi < \phi_1$ there is no epidemic, while if $\phi > \phi_2$ an epidemic may occur. In view of the monotonicity of the model with respect to the parameter $\phi$, this will imply the existence of a $\phi_c \in [\phi_1, \phi_2]$ such that Theorem 3 holds.

To prove the existence of $\phi_1 > 0$ small, we let $L$ be a large integer to be fixed later and, for any site $(z,n) \in \mathcal{G}$, introduce the space–time regions

$$A_{z,n} = (Lz, Ln) + \{[-2L, 2L]^d \times [0, 2L]\}$$

and

$$B_{z,n} = (Lz, Ln) + \{[-L, L]^2 \times [L, 2L]\}.$$

A site $(z, n) \in \mathcal{G}$ is said to be good if $\xi_t(x) = 0$ for any $(x, t) \in B_{z,n}$, that is, there is no infected individual in $B_{z,n}$. When $\lambda < \lambda_c$ and $\phi = 0$, the IRP reduces to a subcritical contact process, which implies, by Section 5 of [12], the existence of an event $F_{z,n}$ measurable with respect to the graphical representation of $\xi_t$ in $A_{z,n}$ such that on $F_{z,n}$ site $(z, n)$ is good. Moreover, for any $\varepsilon > 0$ there exists a large enough $L$, fixed from now on, such that

$$P_\xi((z,n) \text{ is good}) \geq P_\xi(F_{z,n}) \geq 1 - \varepsilon/2 \qquad \text{if } \phi = 0.$$

Now that $L$ is fixed, we can find $\phi_1 > 0$ small such that there is no occurrence of Poisson processes with rate $\phi$ inside $A_{z,n}$ with probability at least $1 - \varepsilon/2$ provided $\phi \leq \phi_1$. In conclusion,

$$P_\xi((z,n) \text{ is good}) \geq 1 - \varepsilon \qquad \text{if } \phi \leq \phi_1.$$

To complete the comparison with oriented percolation, we now position oriented edges between sites in $\mathcal{G}$ as follows: For $(x, n), (y, m) \in \mathcal{G}$ we draw an oriented edge from $(x, n)$ to $(y, m)$ if and only if $n \leq m$ and $A_{x,n} \cap A_{y,m} \neq \varnothing$. In that way, we define a 2-dependent oriented percolation process on $\mathcal{G}$ for which sites are open with probability $1 - \varepsilon$. Since $\varepsilon > 0$ can be chosen arbitrarily small, we can conclude, as in Section 2, that there is no epidemic when $\phi \leq \phi_1$.

To get the existence of $\phi_2 < \infty$, we now deal with the limiting case $\phi = \infty$ for which our process becomes a Richardson model $\xi_t : \mathbb{Z}^d \longrightarrow \{0, \kappa\}$ with growth rate $\lambda \kappa$ (see [11]). If we think of the process as a set-valued process in which the state at time $t$ is the set of infected sites, and start the evolution with one infected cluster at the origin, then the following shape theorem holds ([4], Theorem 6, Chapter 1): there exists a convex set $A$ such that, for any $\delta > 0$, there are constants $C_1 < \infty$ and $\gamma_1 > 0$ such that

(5) $$P_\xi((1 - \delta)tA \subset \xi_t \subset (1 + \delta)tA) \geq 1 - C_1 \exp(-\gamma_1 t).$$

Let $L$ and $\Gamma$ denote two large integers to be fixed later, set

$$\mathcal{H} = \{(z, n) \in \mathbb{Z}^2 : z + n \text{ is even and } n \geq 0\},$$

and, for $z \in \mathbb{Z}$, introduce the spatial box

(6) $$B_z = Lze_1 + [-L, L]^d,$$

where $e_1$ denotes the first unit vector of the $d$-dimensional lattice. We say that a site $(z, n) \in \mathcal{H}$ is occupied if each of the clusters in the spatial box $B_z$ has $\kappa$ infected individuals at time $n\Gamma L$. Then, the shape theorem (5) implies that, for any $\varepsilon > 0$, we can pick $L$ and $\Gamma$ sufficiently large so that the set of occupied sites for the Richardson model $\xi_t$ dominates the set of wet sites in a 1-dependent oriented percolation process with parameter $1 - \varepsilon/2$.

To extend the result to $\phi < \infty$ sufficiently large, we require the following good event, denoted by $G_{z,n}$. For any $x \in B_{z-1} \cup B_{z+1}$, we follow the line $\{x\} \times [n\Gamma L, (n + 1)\Gamma L]$ by going forward in time. Each time we encounter a $\times$, we require the next $\bullet$ to appear before a new $\times$ or $\kappa$-arrow is encountered. This good event assures us that the model with parameters $(\lambda, \kappa, \phi)$ exhibits the same behavior as the Richardson model introduced above in the spatial box $B_{z-1} \cup B_{z+1}$ from time $n\Gamma L$ to time $(n + 1)\Gamma L$. To estimate the probability of our good event, we first observe that, in view of the parameters of the Poisson processes involved in our graphical representation, the probability of going through a $\bullet$ before any $\times$'s or $\kappa$-arrow is given by

$$\frac{(\kappa - 1)\phi}{(\kappa - 1)\phi + 2d\lambda + \kappa}.$$

Let $M$ denote the number of $\times$'s contained in $\{B_{z-1} \cup B_{z+1}\} \times [n\Gamma L, (n + 1)\Gamma L]$. In view of the properties of the exponential distribution, we get

$$E_\xi(M) = \kappa(4L + 1)(2L + 1)^{d-1}\Gamma L$$

since there are $(4L+1)(2L+1)^{d-1}$ sites in $B_{z-1} \cup B_{z+1}$. Finally, by decomposing the event to be estimated according to whether $M > 2E_\xi(M)$ or not, large deviation estimates imply that

$$P_\xi(G_{z,n}) \geq 1 - C_2 \exp(-\gamma_2 \Gamma L)$$
$$- 2\kappa(4L+1)(2L+1)^{d-1} \Gamma L \times \frac{2d\lambda + \kappa}{(\kappa-1)\phi + 2d\lambda + \kappa}$$

for appropriate $C_2 < \infty$ and $\gamma_2 > 0$. In particular, by setting $\phi = \exp(\Gamma L)$, and then choosing $L$ sufficiently large, $P_\xi(G_{z,n})$ can be made greater than $1 - \varepsilon/2$. Putting things together, it follows that the set of occupied sites for the model with parameters $(\lambda, \kappa, \phi)$ dominates the set of wet sites in a 1-dependent oriented percolation process with parameter $1 - \varepsilon$.

To construct our nontrivial stationary distribution, we start $\xi_t$ from the "all $\kappa$" configuration, run the process to time $S$, take the Cesaro average of the distribution at times $0 \leq s \leq S$ and extract a convergent subsequence. By Liggett ([9], Proposition 1.8) the limit $\mu$ is a stationary measure. Since percolation occurs with positive probability when $\varepsilon > 0$ is small enough,

$$\liminf_{n \to \infty} \mu((z,n) \text{ is occupied}) \geq \liminf_{n \to \infty} P_\varepsilon((z,n) \text{ is wet}) > 0,$$

where $P_\varepsilon$ denotes the law of the oriented percolation process with parameter $1 - \varepsilon$, which implies that $\mu$ concentrates on configurations with infinitely many infected clusters.

**4. Proof of Theorem 5.** In this section we prove that if $\lambda \in (0, \lambda_c)$ and $\phi > 1$ there exists a critical value $\kappa_c \geq 2$ such that if $\kappa < \kappa_c$ the IRP converges to the "all 0" configuration, while if $\kappa > \kappa_c$ there is a nontrivial equilibrium. First of all, if $\kappa = 1$ the value of $\phi$ is irrelevant and the process reduces to a basic contact process with parameter $\lambda < \lambda_c$ so that the "all 0" configuration is the unique stationary distribution. In particular, in view of the monotonicity of the IRP with respect to $\kappa$, it suffices to prove the existence of a large $\kappa_0$ such that the process with parameters $(\lambda, \kappa_0, \phi)$ has a nontrivial stationary distribution. This will imply Theorem 5.

To prove the existence of such a $\kappa_0$, the strategy is to compare the particle system with a 1-dependent oriented percolation process on $\mathcal{H}$. To rescale the particle system, we let $T_\kappa = \kappa^2$, and say that a site $(z, n) \in \mathcal{H}$ is occupied if at any time $t \in [nT_\kappa, (n+1)T_\kappa]$ there is at least one infected individual at site $ze_1$, where $e_1$ denotes the first unit vector of the $d$-dimensional lattice. Theorem 5 then follows from the following.

LEMMA 4.1. *Let $\lambda > 0$ and $\phi > 1$. Then, for any $\varepsilon > 0$, there exists a large enough $\kappa$ such that the set of occupied sites dominates the set of wet sites in a 1-dependent oriented percolation process with parameter $1 - \varepsilon$.*

The first step is to investigate the behavior of an isolated cluster when the outside infection rate is set to 0. More precisely, we need some estimates on the first recovery time of a large cluster, the main ingredient (given in Lemma 4.4) being that, with probability close to 1, the infection dies out quickly or persists at least $2T_\kappa$ units of time.

*Random walk estimates.* To figure out the behavior of an isolated cluster, we first consider the continuous-time random walk $\mathcal{X}_t \in \{0, 1, \ldots, \kappa\}$ that makes transitions

$$i \to \begin{cases} i+1, & \text{at rate } i\phi, \\ i-1, & \text{at rate } i, \end{cases} \quad \text{when } 1 \leq i \leq \kappa - 1,$$

and $\kappa \to \kappa - 1$ at rate $\kappa$. In particular, $\mathcal{X}_t$ is equal in distribution to $\xi_t(x)$, the number of infected individuals at site $x \in \mathbb{Z}^d$, when $\lambda = 0$. To estimate the recovery time of a given cluster, we also introduce the stopping time

$$\tau_0 = \inf\{t \geq 0 : \mathcal{X}_t = 0\}$$

and the good event

$$\Omega_\kappa = \{\mathcal{X}_t = \kappa \text{ for some } t < 2T_\kappa\}.$$

The aim is to prove that, with probability close to 1 when $\kappa$ is large, $\tau_0 > 2T_\kappa$ on the event $\Omega_\kappa$ while the random walk returns quickly to 0 on the event $\Omega_\kappa^c$. By convention, all through this section the subscripts and superscripts on the probabilities denote the process and its initial state respectively (e.g., $P_\mathcal{X}^\kappa$ for the law of the process $\mathcal{X}_t$ starting from $\mathcal{X}_0 = \kappa$).

LEMMA 4.2. *Let $\phi > 1$. Then there exist $C_3 < \infty$ and $\gamma_3 > 0$ such that*

$$P_\mathcal{X}^\kappa(\tau_0 < 2T_\kappa) \leq 4\kappa\phi T_\kappa \phi^{-\kappa} + C_3 \exp(-\gamma_3 T_\kappa).$$

PROOF. The basic idea is to deduce the result from similar estimates for the asymmetric discrete-time random walk $X_n$ that makes transitions

$$i \to \begin{cases} i+1, & \text{with probability } \phi/(1+\phi), \\ i-1, & \text{with probability } 1/(1+\phi), \end{cases} \quad \text{when } i \leq \kappa - 1,$$

and $\kappa \to \kappa - 1$ with probability 1. Let $Y_n \in \mathbb{Z}$ be the backward random walk with transitions

$$i \to \begin{cases} i+1, & \text{with probability } 1/(1+\phi), \\ i-1, & \text{with probability } \phi/(1+\phi), \end{cases}$$

and, for any $i \leq 0$, $p_i$ the probability that $Y_n = 0$ for some $n \geq 0$ when starting from $Y_0 = i$. Decomposing according to whether $Y_n$ first jumps to $i - 1$ or $i + 1$ leads to

$$p_i = \frac{\phi}{1+\phi} p_{i-1} + \frac{1}{1+\phi} p_{i+1} \quad \text{and} \quad p_0 = 1.$$

Since $p_i \to 0$ as $i \to -\infty$, it follows that $p_i = \phi^i$, which implies that, for any $n \geq 0$,

$$P_X^\kappa(X_n \leq 0) \leq P_Y^{-\kappa}(Y_j = 0 \text{ for some } j \geq 0) = p_{-\kappa} = \phi^{-\kappa}.$$

Let $K_t$ denote the number of times the process $\mathcal{X}_s$ jumps by time $t$. Since $\mathcal{X}_s$ makes transitions at rate at most $\kappa\phi$, large deviation estimates imply that

$$P_{\mathcal{X}}^\kappa(K_t > 2\kappa\phi t) \leq C_3 \exp(-\gamma_3 t/2)$$

for appropriate $C_3 < \infty$ and $\gamma_3 > 0$. Putting things together and decomposing according to whether $K_{2T_\kappa}$ is smaller or greater than $4\kappa\phi T_\kappa$, we obtain

$$P_{\mathcal{X}}^\kappa(\tau_0 < 2T_\kappa) \leq P_{\mathcal{X}}^\kappa(\tau_0 < 2T_\kappa; K_{2T_\kappa} < 4\kappa\phi T_\kappa) + P_{\mathcal{X}}^\kappa(K_{2T_\kappa} \geq 4\kappa\phi T_\kappa)$$
$$\leq P_X^\kappa(X_n \leq 0 \text{ for some } n < 4\kappa\phi T_\kappa) + P_{\mathcal{X}}^\kappa(K_{2T_\kappa} \geq 4\kappa\phi T_\kappa)$$
$$\leq \sum_{n=1}^{4\kappa\phi T_\kappa} P_X^\kappa(X_n \leq 0) + C_3 \exp(-\gamma_3 T_\kappa)$$
$$\leq 4\kappa\phi T_\kappa \phi^{-\kappa} + C_3 \exp(-\gamma_3 T_\kappa).$$

This completes the proof of the lemma. $\square$

LEMMA 4.3. *Let $\phi > 1$. Then there exist $C_4 < \infty$ and $\gamma_4 > 0$ such that for any $t < 2T_\kappa$*

$$P_{\mathcal{X}}^1(\tau_0 > t; \Omega_\kappa^c) \leq C_4 \exp(-\gamma_4 t).$$

PROOF. To begin with, we prove the analogous result for the asymmetric random walk $Z_n \in \mathbb{Z}$ that makes transitions

$$i \to \begin{cases} i+1, & \text{with probability } \phi/(1+\phi), \\ i-1, & \text{with probability } 1/(1+\phi). \end{cases}$$

Let $\sigma_0$ denote the first time $Z_n = 0$ and $p = \phi/(1+\phi)$. First of all, since

$$E_Z^1[Z_n - Z_0] = (2p-1)n,$$

large deviation estimates on $Z_n$ imply that, for any $\delta > 0$,

$$P_Z^1(Z_n \leq (2p-1-\delta)n) \leq C_5 \exp(-\gamma_5 n)$$

for appropriate $C_5 < \infty$ and $\gamma_5 > 0$. In other respects, since $\phi > 1$, and so $p > 0.5$, we can find a sufficiently small $\delta > 0$ such that $2p - 1 - \delta > 0$, which implies that

$$P_Z^1(n < \sigma_0 < \infty) = \sum_{k=n+1}^\infty P_Z^1(\sigma_0 = k)$$

$$\leq \sum_{k=n+1}^{\infty} P_Z^1(Z_k \leq 0)$$

$$\leq C_5 \sum_{k=n+1}^{\infty} \exp(-\gamma_5 k)$$

$$\leq C_6 \exp(-\gamma_5 n).$$

To deduce the result for $\mathcal{X}_s$, we first observe that, on $\Omega_\kappa^c$ and by time $\min(\tau_0, 2T_\kappa)$, the continuous-time random walk $\mathcal{X}_s$ is just a time change of $Z_n$. Since $\mathcal{X}_s$ jumps at rate at least 1 by time $\tau_0$, the probability that $K_t < t$ (where $K_t$ denotes the number of jumps by time $t$) when $\tau_0 > t$ is smaller than $C_7 \exp(-\gamma_7 t)$. In conclusion, for any $t < 2T_\kappa$, we get

$$P_\mathcal{X}^1(\tau_0 > t; \Omega_\kappa^c) \leq P_\mathcal{X}^1(\tau_0 > t; \Omega_\kappa^c; K_t \geq t) + P_\mathcal{X}^1(\Omega_\kappa^c; K_t < t)$$
$$\leq P_Z^1(t < \sigma_0 < \infty) + P_\mathcal{X}^1(\Omega_\kappa^c; K_t < t)$$
$$\leq C_6 \exp(-\gamma_5 t) + C_7 \exp(-\gamma_7 t) \leq C_4 \exp(-\gamma_4 t)$$

for appropriate $C_4 < \infty$ and $\gamma_4 > 0$. This completes the proof. □

LEMMA 4.4. *Let $\phi > 1$ and $T_\kappa = \kappa^2$. Then there exist $C_8 < \infty$ and $\gamma_8 > 0$ such that*

$$P_\mathcal{X}^1(t < \tau_0 < 2T_\kappa) \leq C_8 \exp(-\gamma_8 \kappa) + C_8 \exp(-\gamma_8 t).$$

PROOF. To begin with, we observe that, on the event $\Omega_\kappa$, the first time $\mathcal{X}_s$ hits 0 when $\mathcal{X}_0 = \kappa$ is bounded in distribution by the first time it hits 0 when $\mathcal{X}_0 = 1$ so that

$$P_\mathcal{X}^1(\tau_0 < 2T_\kappa; \Omega_\kappa) \leq P_\mathcal{X}^\kappa(\tau_0 < 2T_\kappa).$$

In particular, by decomposing the event to be estimated according to whether $\Omega_\kappa$ occurs or not, it follows from Lemma 4.2 and Lemma 4.3 that

$$P_\mathcal{X}^1(t < \tau_0 < 2T_\kappa) \leq P_\mathcal{X}^1(t < \tau_0 < 2T_\kappa; \Omega_\kappa) + P_\mathcal{X}^1(t < \tau_0 < 2T_\kappa; \Omega_\kappa^c)$$
$$\leq P_\mathcal{X}^1(\tau_0 < 2T_\kappa; \Omega_\kappa) + P_\mathcal{X}^1(\tau_0 > t; \Omega_\kappa^c)$$
$$\leq 4\kappa^3 \phi^{1-\kappa} + C_3 \exp(-\gamma_3 \kappa^2) + C_4 \exp(-\gamma_4 t)$$
$$\leq C_8 \exp(-\gamma_8 \kappa) + C_8 \exp(-\gamma_8 t)$$

for suitable $C_8 < \infty$ and $\gamma_8 > 0$. This completes the proof of the lemma. □

PROOF OF LEMMA 4.1. To deduce Lemma 4.1, we show that if the cluster at site 0 has at least one infected individual from time 0 to time $T_\kappa$, that is, site $(0,0)$ is occupied, then this individual gives birth at site $z \sim 0$ and by time $T_\kappa$ to an infection that survives at least $2T_\kappa$ units of time, so that $(z, 1)$ is occupied. The probability of this event has to be greater than $1 - \varepsilon$ for $\kappa$ sufficiently large.

Each time the infected individual at site 0 gives rise to a new infection at site $z$, we call this infection a strong infection if it lives at least $2T_\kappa$ units of time (at site $z$). Let $\sigma^z$ denote the first time a strong infection is born at site $z$, that is,

$$\sigma^z = \inf\{t \geq 0 : \xi_s(z) \neq 0 \text{ for any } t \leq s \leq t + 2T_\kappa\}.$$

Then, Lemma 4.1 follows from the following.

LEMMA 4.5. *Let $\lambda > 0$, $\phi > 1$ and $T_\kappa = \kappa^2$. Then $P_\xi(\sigma^z > T_\kappa) \to 0$ as $\kappa \to \infty$.*

PROOF. The proof relies on a restart argument. The basic idea is that each time the cluster at site $z$ gets infected (between time 0 and time $T_\kappa$), the infection dies out quickly or lives at least $2T_\kappa$ units of time. In particular, the number of trials before having a strong infection can be made arbitrarily large by choosing $\kappa$ large enough. Since a geometrical number of trials suffices, the first strong infection will appear by time $T_\kappa$ with probability close to 1. To make this argument precise, we introduce, for any $i \geq 1$, the stopping times

$$\rho_i = \inf\{t \geq \bar{\rho}_{i-1} : \xi_t(z) = 1\} \quad \text{and} \quad \bar{\rho}_i = \inf\{t \geq \rho_i : \xi_t(z) = 0\}$$

with the convention $\bar{\rho}_0 = 0$. That is, $\rho_i$ is the $i$th time the cluster at site $z$ gets infected, and $\bar{\rho}_i$ is the $i$th time the cluster recovers. Let

$$N = \inf\{i \geq 1 : \bar{\rho}_i - \rho_i > 2T_\kappa\}$$

so that $\rho_N = \sigma^z$, the first time a strong infection is born at site $z$. By decomposing the event to be estimated according to the number of trials $N$, we get

$$P_\xi(\sigma^z > T_\kappa) = \sum_{n=1}^{\infty} P_\xi(\rho_n > T_\kappa | N = n) P_\xi(N = n)$$

$$\leq \sum_{n=1}^{\infty} P_\xi(\rho_{i+1} - \rho_i > T_\kappa/n \text{ and } \bar{\rho}_i - \rho_i \leq 2T_\kappa \text{ for some } i \leq n-1)$$

$$\times P_\xi(N = n)$$

$$\leq \sum_{n=1}^{\infty} \sum_{i=1}^{n-1} P_\xi(\rho_{i+1} - \rho_i > T_\kappa/n \text{ and } \bar{\rho}_i - \rho_i \leq 2T_\kappa) P_\xi(N = n)$$

$$\leq \sum_{n=1}^{\infty} n P_\xi(\rho_2 - \rho_1 > T_\kappa/n \text{ and } \bar{\rho}_1 - \rho_1 \leq 2T_\kappa) P_\xi(N = n).$$

To prove that the previous series can be made arbitrarily small by taking $\kappa$ large enough, we first observe that, since $\mathcal{X}_t$ has a drift to the right, there is a constant $q = q(\phi) > 0$ such that

$$P^1_{\mathcal{X}}(\Omega_\kappa) > q \qquad \text{for any } \kappa \geq 1.$$

This, together with Lemma 4.2, implies that $P^1_{\mathcal{X}}(\tau_0 > 2T_\kappa) = q_\kappa > q/2$ for all $\kappa$ sufficiently large. In particular, for any $\varepsilon > 0$, there exists a large enough $n_\varepsilon$ such that

$$\sum_{n=n_\varepsilon}^\infty n P_\xi(N = n) = \sum_{n=n_\varepsilon}^\infty n q_\kappa (1 - q_\kappa)^{n-1} \leq \varepsilon/2.$$

To show that the first $n_\varepsilon$ terms tend to 0 as $\kappa \to \infty$, we first observe that

$$P_\xi(\rho_2 - \rho_1 > T_\kappa/n \text{ and } \bar\rho_1 - \rho_1 \leq 2T_\kappa)$$
$$\leq P_\xi(\rho_2 - \bar\rho_1 > T_\kappa/2n) + P_\xi(T_\kappa/2n < \bar\rho_1 - \rho_1 \leq 2T_\kappa).$$

Since there is at least one infected individual at site 0, we get

$$P_\xi(\rho_2 - \bar\rho_1 > T_\kappa/2n) \leq \exp(-\lambda T_\kappa/2n_\varepsilon) \qquad \text{for } n \leq n_\varepsilon.$$

In other respects, Lemma 4.4 implies that $P_\xi(T_\kappa/2n < \bar\rho_1 - \rho_1 \leq 2T_\kappa)$ tends to 0 as $\kappa \to \infty$. In conclusion,

$$\sum_{n=1}^{n_\varepsilon} n P_\xi(\rho_2 - \rho_1 > T_\kappa/n \text{ and } \bar\rho_1 - \rho_1 \leq 2T_\kappa) \leq \varepsilon/2$$

for $\kappa$ sufficiently large, which completes the proof. $\square$

Lemma 4.5 implies that if there is an infected individual at site 0 at any time $t \in [0, T_\kappa]$, then, with probability close to 1 when $\kappa$ is large, the clusters at sites $z$, $z \sim 0$, have each at least one infected individual at any time $t \in [T_\kappa, 2T_\kappa]$, which proves Lemma 4.1. The existence of a nontrivial stationary distribution can then be deduced from Lemma 4.1 as in Section 3. $\square$

**5. Proof of Theorem 6.** In this section we now deal with the CRP, that is, the recovery mechanism is now described by the transitions $i \to 0$, $i = 1, 2, \ldots, \kappa$, at rate 1. We prove that, contrary to the IRP, for any within infection rate $\phi \geq 0$, we can find a critical value $\lambda_c(\phi)$ such that, provided $\lambda \leq \lambda_c(\phi)$, there is no epidemic whatever the cluster size. As explained in Section 1, the key idea of the proof is that the state of each cluster returns to 0 sufficiently often so that, between two consecutive recoveries, the number of infected individuals in a given cluster cannot exceed some threshold $\kappa_0$ with probability close to 1.

As previously, we rely on a rescaling argument to compare the particle system with the oriented percolation process introduced in Section 1. For any $(z, n) \in \mathcal{G} = \mathbb{Z}^d \times \mathbb{Z}_+$, let $S_{z,n}$ be the vertical segment and let $C_{z,n}$ be the space–time cylinder given by (2) and (3), respectively. The good event, denoted by $H_{z,n}$, we now consider has to be measurable with respect to the graphical representation in $C_{z,n}$ and assure us that, at time $(n+1)T$, the cluster at site $z$ is healthy.

To construct the good event $H_{z,0}$, we first require the clusters at sites $x$, $x \sim z$, to recover at least once between time 0 and time $T/2$ (which corresponds to at least one occurrence of Poisson processes with rate 1 by time $T/2$), and the cluster at site $z$ to recover at least once between time $T/2$ and time $T$. The probability of this event can be estimated

$P_\eta($there is at least one $\times_1$ on the segment $\{x\} \times [0, T/2]$ for any $x \sim z$

and at least one $\times_1$ on the segment $\{z\} \times [T/2, T])$

$\geq 1 - (2d+1) \exp(-T/2)$,

where $P_\eta$ denotes the law of the CRP. In particular, given $\varepsilon > 0$, there is a large enough $T$, fixed from now on, so that the previous event has probability at least $1 - \varepsilon/3$. The reason why we want the neighboring clusters of $z$ to recover at least once by time $T/2$ is to control the number of infected individuals around $z$ until time $T$; the aim is to prevent infections at site $z$ coming from the outside.

To estimate the number of infected individuals in the neighborhood of $z$ until time $T$, we introduce the continuous-time random walk $\mathcal{Z}_t \in \mathbb{Z}_+^*$ that makes transitions

$$i \to \begin{cases} i+1, & \text{at rate } i\phi, \\ 1, & \text{at rate } 1. \end{cases}$$

By monotonicity of the sequence $\{\mathcal{Z}_t \leq n, \forall\, t \leq T\}$, $n \geq 1$, the Beppo–Lévi theorem implies that

$$\lim_{n \to \infty} P_\mathcal{Z}^1(\mathcal{Z}_t \leq n \text{ for all } t \leq T) = P_\mathcal{Z}^1(\mathcal{Z}_t < \infty \text{ for all } t \leq T) = 1,$$

where $P_\mathcal{Z}^1$ denotes the law of the process $\mathcal{Z}_t$ starting from $\mathcal{Z}_0 = 1$. In particular, we can find a sufficiently large $\kappa_0 < \infty$ depending on $\phi$ and $T$, fixed from now on, such that

$$P_\mathcal{Z}^1(\mathcal{Z}_t \geq \kappa_0 \text{ for some } 0 \leq t \leq T) \leq \varepsilon/6d.$$

Moreover, since the neighboring clusters of $z$ recover at least once by time $T/2$, for any $x \sim z$ and any time $t \in [T/2, T]$, we get $\eta_t(x) \leq \mathcal{Z}_t$ in distribution so that

$$P_\eta(\eta_t(x) \geq \kappa_0 \text{ for some } x \sim z \text{ and some } T/2 \leq t \leq T) \leq \varepsilon/3.$$

We now fix $\lambda_c(\phi)$ so that, with high probability, no epidemic may occur provided $\lambda \leq \lambda_c(\phi)$. When $\kappa \leq \kappa_0$, the result follows from Theorem 2 by taking $\lambda_c(\phi) = \lambda_c/\kappa_0$. To deal with the nontrivial case $\kappa > \kappa_0$, we complete the construction of our good event $H_{z,0}$ by requiring no infection of site $z$ coming from the outside, that is, from a site $x$ with $x \sim z$, until time $T$. Since each of the neighboring clusters of $z$ has at most $\kappa_0$ infected individuals, this occurs if no $i$-arrow with $1 \leq i \leq \kappa_0$ points at site $z$ between time $T/2$ and time $T$, an event with probability

$$\exp(-d\lambda\kappa_0 T),$$

which can be made greater than $1 - \varepsilon/3$ by choosing $\lambda > 0$ sufficiently small. Since there is no outside infection at site $z$ and the cluster at $z$ recovers at least once between times $T/2$ and $T$, the event $H_{z,0}$ assures us that $z$ is healthy at time $T$. Finally, the probability of $H_{z,0}$ being greater than $1 - \varepsilon$ with $\varepsilon > 0$ arbitrarily small, Theorem 6 follows from a comparison with the oriented percolation process introduced in Section 2.

**Acknowledgments.** The authors wish to thank Ellen Saada and Claudia Neuhauser for the time they devoted to them and fruitful discussions, and an anonymous referee for his suggestions to simplify some of the proofs.

Laboratoire de Mathématiques Raphaël Salem
UMR 6085, CNRS—Université de Rouen
Avenue de l'Université
BP. 12
76801 Saint Etienne du Rouvray
France
e-mail: Lamia.Belhadji@univ-rouen.fr

Department of Ecology
 Evolution, and Behavior
University of Minnesota
Ecology Building
1987 Upper Buford Circle
St. Paul, Minnesota 55108
USA
e-mail: Nicolas.Lanchier@univ-rouen.fr